\documentclass[12pt]{article}
\usepackage[cp1251]{inputenc}
\usepackage[russian]{babel}
\usepackage{amsfonts}

\newtheorem{lemma}{Lemma}
\newtheorem{theorem}{Theorem}

\begin{document}

\title{On Hamiltonian-Minimal Lagrangian Tori in~${\Bbb C}P^2$}
\author{A.E. Mironov}
\date{}
\maketitle

\section{Introduction}

In this article we obtain some equations for
Hamiltonian-minimal Lagrangian surfaces in~${\Bbb CP}^2$ (Lem\-ma~1)
and give their particular solutions in the case of tori.

An~immersion $\psi: \Omega\rightarrow {\Bbb CP}^2$
of a~domain $\Omega\subset{\Bbb R}^2$
is {\it Lagrangian\/} if $\psi^*(S)=0$
where $S$
is the Fubini--Study form of~${\Bbb CP}^2$.
An~immersion is $H$-{\it minimal\/} if the area variations
along all Hamiltonian fields are zero.
A~field $W$ is Hamiltonian along $\psi(\Omega)$
if the 1-form $S(W,.)$ is exact on~$\psi(\Omega)$.
Our main result is the following

\begin{theorem} The mapping $\psi:{\Bbb R}^2\rightarrow {\Bbb
CP}^2$ defined by the formula
$$
\psi(x,y)=(F_1(x)e^{i(G_1(x)+\alpha_1y)}: F_2(x)e^{i(G_2(x)+\alpha_2y)}:
F_3(x)e^{i(G_3(x)+\alpha_3y)}),
$$
is a~conformal Lagrangian $H$-minimal immersion of the plane.
Here
$$
 F_i=\sqrt{\frac{e^{2v(x)}+\alpha_{i+1}\alpha_{i+2}}
 {(\alpha_i-\alpha_{i+1})(\alpha_i-\alpha_{i+2})}},
\quad
 G_i=\frac{\alpha_i}{2}\int\limits_{x_0}^{x}
\frac{2c_2-ae^{2v(z)}}{\alpha_ie^{2v(z)}-c_1}dz,
$$
$$
 e^{2v(x)}=a_1\left(1-\frac{a_1-a_2}{a_1}
\rm{sn}^2\left(x\sqrt{a_1+a_3},\frac{a_1-a_2}
 {a_1+a_3}\right)\right)
 \eqno{(1)}
$$
$($the index $i$ is considered modulo~$3)$, $a_1>a_2>0$,
$\alpha_i$ are reals satisfying ~{\rm (15)} and~{\rm(16)}, the
constants~$c_1$, $c_2$, $a$, and~$a_3$ are expressed in terms
of~$a_i$ and~$\alpha_i$ by {\rm (10), (12)}, and {\rm~(13)}, and
$\rm{sn}$ is the Jacobi elliptic function.

If, in addition,
$\alpha_1,\alpha_2,\alpha_3\in {\Bbb Z }$;
$\lambda_1=G_1(T)-G_3(T)+(\alpha_1-\alpha_3)\tau$,
$\lambda_2=G_2(T)-G_3(T)+
(\alpha_2-\alpha_3)\tau\in{2 \pi \Bbb Q}$,
where $T$ is the period of the function~$e^{2v(x)}$
$($see~$(14))$,
and
$\tau\in{\Bbb R}$;
then $\psi$ is a~doubly periodic mapping with periods
$e_1=(0,1)$ and $e_2=N(T,\tau)$,
where
$N$
is some natural number.
\end{theorem}

Observe that $\lambda_1$ and $\lambda_2$ depend on
the free parameters $a_1$, $a_2$, and~$\tau$; therefore,
$\lambda_1,\lambda_2\in{2 \pi \Bbb Q}$
for a~dense set of triples~$(a_1,a_2,\tau)$ in some domain.

The notion of $H$-minimality was introduced in~[1];
in the same article it was proven that the Clifford tori in~
${\Bbb C}^n$ are Lagrangian $H$-minimal.
Other examples of such tori in~${\Bbb C}^2$ are constructed in~[2,\,3].
In~[4] minimal Lagrangian tori in~${\Bbb CP}^2$ are constructed
which are particular cases of the tori of Theorem~1
for $\alpha_1+\alpha_2+\alpha_3=0$ and
$
 (a_1+a_2)\bigl(c_1^2+c_2^2\bigr)-a^2_1a^2_2=0.
$

The author is grateful to M.~V. Neshchadim
who brought system~(3)--(5) to a~lucid form
which made it possible to derive from ~(3)--(5)
the 
 Tzitz\'eica
equation for the metric of a~minimal Lagrangian torus~in~
${\Bbb C}P^2$.

\section{Proof of Theorem~1}

Let $S^5$ be the unit sphere in~${\Bbb C}^3$
and let
${\cal H}:S^5\rightarrow {\Bbb CP}^2$
be the Hopf bundle.
Denote by $\omega$ the symplectic form on~${\Bbb C}^3$:
$$
 \omega=dx_1\wedge dy_1+dx_2\wedge dy_2+dx_3\wedge dy_3,
$$
where
$z_j=x_j+iy_j$
are coordinates in~${\Bbb C}^3$, $j=1,2,3$.
Suppose that $L$ is a~surface in~${\Bbb CP}^2$
and
${\cal U}$
is a~sufficiently small neighborhood of a~point $p\in L$. Denote by
$\widetilde{{\cal U}}$
some horizontal lift of ${\cal U}$ to~$S^5$.

A~criterion for $L$ to be Lagrangian is as follows
(see~[5]):
{\sl the surface~$L$ is Lagrangian if and only if the linear span
of the radius vector~$\tilde{p}$
$(\tilde{p}$ is the lift of~$p)$
and the tangent plane to~$\widetilde{{\cal U}}$ at~
$\tilde{p}$ is a~Lagrangian three-dimensional subspace in~
${\Bbb C}^3$ for all~$p\in L$.}

We also use the following criterion for $H$-minimality
in terms of the Lagrangian angle: 
{\sl the surface~$L$
is $H$-minimal if and only if the Lagrangian angle is a~harmonic
function on~$L$ in the induced~ metric.}

The Lagrangian angle is a~function on $L$ constructed (locally) as follows:
Take an~orientation on~$\widetilde{{\cal U}}$.
Put
$$
 e^{-i\beta}=z_1\wedge z_2\wedge z_3(\xi_1,\xi_2,p),
$$
where $\xi_1$ and  ~$\xi_2$
make an~orthonormal tangent  basis for~$\widetilde{{\cal U}}$
agreeing with orientation. The function $\beta(p)$
is called the {\it Lagrangian angle}.
In general, $\beta(p)$
is a~many-valued function on~$L$.
It may change its values by $2\pi k$, $k\in{\Bbb Z}$,
in passing around a~cycle.

We define the conformal Lagrangian immersion~$\psi$ of a~domain
$\Omega\subset{\Bbb R}^2$ with coordinates~$x$ and~$y$
into~${\Bbb CP}^2$ as the composite of
$r:\Omega\rightarrow S^5$ and ${\cal H}$.
Let $e^{2v(x,y)}(dx^2+dy^2)$
be the induced metric on~$\psi(\Omega)$.
Note that, since $\psi$ is conformal and Lagrangian, we have
$$
\langle r,r_x\rangle =\langle r,r_y\rangle =
\langle r_x,r_y\rangle =0,\quad  |r_x|=|r_y|=e^{v},
$$
where $\langle .,.\rangle $
is the Hermitian product.
Consequently,
$$
 R=
  \left(
  \begin{array}{ccc}
   e^{i\beta} r^1 & e^{i\beta} r^2 & e^{i\beta} r^3\\
   e^{-v}r^1_x & e^{-v}r^2_x & e^{-v}r^3_x \\
   e^{-v}r^1_y & e^{-v}r^2_y & e^{-v}r^3_y \\
  \end{array}\right)\in{\rm SU(3)},
$$
where $r^1$, $r^2$, and $r^3$ are the components of~$r$.
Thus, for the matrices~$A$ and $B$,
$$
 A=
  \left(
  \begin{array}{ccc}
   i\beta_x & e^{v+i\beta} & 0\\
  -e^{v-i\beta} & -if-i\beta_x & ig-v_y \\
  0 & ig+v_y & if\\
  \end{array}\right),
$$
$$
 B=
  \left(
  \begin{array}{ccc}
   i\beta_y & 0 & e^{v+i\beta} \\
  0 & ig & if+v_x  \\
  -e^{v-i\beta} & if-v_x & -ig-i\beta_y\\
  \end{array}\right),
$$
in the Lie algebra ${\rm{su}(3)}$,
where
$$if=\langle \partial_x(e^{-v}r_y),e^{-v}r_y\rangle $$
and
$$ig=\langle \partial_y(e^{-v}r_x),e^{-v}r_x\rangle, $$
the following equalities hold:
$$
R_x=AR,\quad  R_y=BR.
                                                    \eqno{(2)}
$$
The matrices $A$ and $B$ satisfy the zero curvature equation
$$
 A_y-B_x+[A,B]=0.
$$
The nonzero components of this equation have the form
$$
 f_y+2fv_y+g_x+2gv_x+\beta_{xy}=0,
$$
$$
-e^{2v}+2f^2+2g^2+ig_y
+iv_y\beta_y+g(2iv_y+\beta_y)
-v_{yy}-if_x-iv_x\beta_x
$$
$$
+f(-2iv_x+\beta_x)-v_{xx}=0,
$$
$$
e^{2v}-2f^2-2g^2+ig_y+iv_y\beta_y+g(2iv_y-\beta_y)
+v_{yy}-if_x-iv_x\beta_x
$$
$$
+f(-2iv_x-\beta_x)+v_{xx}=0.
$$
Hence, we obtain

\begin{lemma}
The equations hold:
$$
U_y+V_x+e^{2v}\beta_{xy}=0,\eqno{(3)}
$$
$$
V_y+v_ye^{2v}\beta_y=U_x+v_xe^{2v}\beta_x,\eqno{(4)}
$$
$$
\Delta v+e^{2v}-2(U^2+V^2)e^{-4v}-(\beta_xU+\beta_yV)e^{-2v}=0,\eqno{(5)}
$$
where $U=fe^{2v}$ and $V=ge^{2v}$.
\end{lemma}

Since harmonic functions remain harmonic under conformal changes of the metric
and since harmonic functions on a~torus are constant,
we may assume that the Lagrangian angle for the tori
has the form $\beta=ax+by$, $a,b\in{\Bbb R}$.

Below we consider the case in which the functions~$v$, $f$, and $g$
depend only on~$x$. Then from (3)--(5) we obtain
$$
 g=c_1e^{-2v(x)},\quad f=c_2e^{-2v(x)}-\frac{a}{2},
$$
$$
(v')^2=-\frac{a}{4}-\bigl(c_1^2+c_2^2\bigr)e^{-4v}
+(ac_2-bc_1)e^{-2v}-e^{2v}-c,
                                               \eqno{(6)}
$$
where $c$, $c_1$, and~$c_2$ are some constants.
We will seek ~$r^i$ in the form
$$
 r^i=C_i(x)e^{i\alpha_i y},
$$
where $C_i(x)$
is a~complex-valued function and $\alpha_i\in{\Bbb R}$.
From ~(2) we obtain
$$
 2(e^{4v(x)}+c_1\alpha_i)C_i(x)+iC_i'(x)
(2c_2+ae^{2v(x)}+2ie^{2v(x)}v'(x))+2e^{2v}C''_i(x)=0,
                                                           \eqno{(7)}
$$
$$
 2i(c_1-e^{2v(x)}\alpha_i)C'(x)+\alpha_iC(x)
((a+2iv'(x))e^{2v(x)}-2c_2)=0,\eqno{(8)}
$$
$$
 2\bigl(e^{2v(x)}\bigl(e^{2v(x)}-b\alpha_i-\alpha_i^2\bigr)
-c_1\alpha_i\bigr)C_i(x)+
 C'(x)((ia+2v'(x))e^{2v(x)}-2ic_2)=0.
                                                         \eqno{(9)}
$$
Note that if $\alpha_i$ satisfies the equation
$$
 \alpha^3+b\alpha^2+c\alpha+c_1=0
$$
then ~(7) and (9) ensue from ~(8) and (6).
Now, from the condition $R\in{\rm{SU}(3)}$ and~(8) we find
that
$$
 C_i(x)=F_i(x)e^{iG_i(x)},
$$
where
$$
c_1=-\alpha_1\alpha_2\alpha_3,\quad
c=\alpha_1\alpha_2+\alpha_1\alpha_3+\alpha_2\alpha_3,
\quad  b=-\alpha_1-\alpha_2-\alpha_3.
                                                        \eqno{(10)}
$$
It remains to find a~solution to ~(6).
Executing the change of variables $h=e^{2v}$,
we take equation~(6) to the form
$$
 (h')^2+4(h-a_1)(h-a_2)(h+a_3)
$$
$$
 =(h')^2+4h^3+(4c+a^2)h^2+4(bc_1-ac_2)h+4(c_1^2+c_2^2)=0,
  \eqno{(11)}
$$
where
$$
a_3=\frac{c_1^2+c_2^2}{a_1a_2},\quad
a=\frac{bc_1+a_1a_3+a_2a_3-a_1a_2}{c_2},
                                                               \eqno{(12)}
$$
and $c_2$ is a~root of the equation
$$
 c^4_2(a_1-a_2)^2+2c^2_2
\bigl(a_1^3a_2^2+a_1^2a_2^3+\bigl(a_1a_2^2+a_1^2a_2\bigr)bc_1
+\bigl(a_1^2+a_2^2\bigr)c_1^2+2a_1^2a_2^2c\bigr)
$$
$$
+\bigl((a_1+a_2)c_1^2-a_1^2a_2^2+a_1a_2c_1b\bigr)^2=0.
\eqno(13)
$$
From the identity
$$
(\rm{sn}(x)')^2
=(1-\rm{sn}^2(x))(1-k^2\rm{sn}^2(x))
$$
(see~[6]) we  infer easily that ~(11)
has a~solution of the form~(1), with
$$
 \rm{sn}(x,k)=\sin\varphi
$$
and $\varphi$
is the inverse function of~
$$
 w(\varphi)=\int\limits_0^{\varphi}\frac{dt}{\sqrt{1-k^2\sin^2t}},
\quad
 0<k<1.
$$
The function~$e^{2v(x)}$ has the period
$$
T=\frac{2w(\frac{\pi}{2})}{\sqrt{a_1+a_3}}.\eqno{(14)}
$$
The choice of the parameters $a_i$ and~$\alpha_i$
is restricted by the condition $c_2\in{\Bbb R}$:
$$
 P=a_1^3a_2^2+a_1^2a_2^3+
\bigl(a_1a_2^2+a_1^2a_2\bigr)bc_1+\bigl(a_1^2+a_2^2\bigr)c_1^2
+2a_1^2a_2^2c\leq 0,
                                                            \eqno{(15)}
$$
$$
  P^2-(a_1-a_2)^2\bigl((a_1+a_2)c_1^2-a_1^2a_2^2+a_1a_2c_1b\bigr)^2\geq 0.
                                                               \eqno{(16)}
$$
Theorem~1 is proven.

Inequalities~(15) and (16) are satisfied, for example, for
$a_1=2$, $a_2=1$, $\alpha_1=0$, $\alpha_2=-1$, and $\alpha_3=3$.

\begin{center} \bf{References} \end{center}

\vskip3mm

[1] Oh, Y. Volume minimization of Lagrangian submanifolds under
Hamil\-to\-nian deformations // Math. Z. 1993. V. 212. P. 175-192.

[2] Castro, I., and Urbano, F. Examples of unstable Hamiltonian-minimal
Lagrangian tori in ${\mathbb C}^2$ // Compositio Math. 1998. V. 111.
P. 1-14.

[3]  Helein, F., and Romon, P. Hamiltonian stationary Lagrangian
surfaces in ${\mathbb C}^2$ //
Comm. Anal. Geom. 2002. V. 10. P. 79--126.

[4] Castro, I., and Urbano, F. New examples of minimal Lagrangian
tori in the complex projective plane // Manuscripta Math. 1994.
V. 85. P. 265--281.

[5] Helein, F., and Romon, P. Hamiltonian stationary Lagrangian
surfaces in Hermitian symmetric spaces // In: Differential
geometry and Integrable Systems. Eds. M. Guest, R. Miyaoka, and
Y. Ohnita. Contemporary Mathe\-matics. V. 308. Amer. Math. Soc.,
Providence, 2002. P. 161-178.

[6]Akhiezer~N.~I., Elements of the Theory of Elliptic Functions //
Moscow: Nauka, 1970 (in Russian).

\vskip3mm

Andrey Mironov

Sobolev Institute of Mathematics, pr. ac. Koptyuga 4,
630090, Novosibirsk, Russia

E-mail address: mironov@math.nsc.ru

\end{document}